\documentclass[
12pt,reqno]{amsart}
\usepackage{verbatim}
\usepackage{amssymb}

\usepackage{enumerate}
\numberwithin{equation}{section}

\newtheorem{theorem}{Theorem}[section]
\newtheorem{lemma}[theorem]{Lemma}
\newtheorem{corollary}[theorem]{Corollary}

\theoremstyle{definition}
\newtheorem{definition}[theorem]{Definition}

\theoremstyle{remark}

\newcommand{\itemprefix}{}
\makeatletter
\newcommand{\myitem}{%
\item\protected@edef\@currentlabel{\itemprefix\theenumi}%
}
\def\Int{\mathop{\operator@font Int}\nolimits}
\makeatother

\newcommand{\mc}[1]{\mathcal{#1}}

\newcommand{\setm}{\setminus}

\newcommand{\subs}{\subset}
\newcommand{\dom}{\operatorname{dom}}

\newcommand{\supp}{\operatorname{sp}}

\author[I. Juh\'asz]{Istv\'an Juh\'asz}
\address      { Alfr\'ed Rényi Institute of Mathematics%
}
\email{juhasz@renyi.hu}

\author[J. van Mill]{Jan van Mill}
\address{University of Amsterdam}
\email{j.vanMill@uva.nl}

\thanks{The first author was supported
by NKFIH grant no. K 129211.}

\title{On the double density spectra of compact spaces}

\begin{document}

\begin{abstract}
The set $dd(X)$ of densities of {\em all} dense subspaces of a topological space $X$ is
called the {\em double density spectrum} of $X$. In this note we present a couple of
results that imply $\lambda \in dd(X)$, provided that $X$ is a {\em compact} space
and $\lambda$ is a cardinal satisfying certain conditions.

As a consequence of these results, we prove that $$dd(X) = [d(X), w(X)]$$ holds for any polyadic space $X$.
This, in turn, implies that $dd(G) = [d(G), w(G)]$ for any
locally compact topological group $G$.
\end{abstract}

\maketitle

\section{Introduction}

For any topological space $X$ we denote by $\mathcal{D}(X)$ the family of all dense subspaces of $X$.
Then $$dd(X) = \{d(D) : D \in \mathcal{D}(X)\},$$ i.e.
the set of densities of {\em all} dense subspaces of $X$ is
called the {\em double density spectrum} of $X$. This concept was introduced and studied in \cite{JvMSSz}.
While the double density spectra of Hausdorff, resp. regular spaces were completely characterised there,
only some partial results were obtained about the double density spectra of compact spaces, in fact only
for separable compacta. In this note we shall present some new results on this problem,
by providing information about the double density spectra of certain well-studied classes of compact spaces,
for instance polyadic spaces.

Of course, we have $d(X) = \min dd(X)$ and $\delta(X) = \sup dd(X) \le \pi(X)$ for any space $X$.
Now, by the main result of \cite{JSh} we have $\pi(X) = \max dd(X)$ if $X$ is compact,
so this problem reduces to the following question: For what cardinals $\lambda$ does $d(X) < \lambda < \pi(X)$
imply $\lambda \in dd(X)$ for a compactum $X$?

For all undefined notions, see \cite{juhasz}.

\section{The class $\Gamma(\lambda)$}

We start by presenting a relatively simple lemma that will turn out to be very useful below.

\begin{lemma}\label{lm:EF}
Assume that $f : X \to Y$ is continuous and $E \subs f[X]$ is such that $|E| = d(E) = \lambda$,
moreover $d(f^{-1}(y)) \le \lambda$ for every $y \in E$. Then there is a dense subset
$F$ of $f^{-1}[E]$ such that $|F| = d(F) = \lambda$ as well.
\end{lemma}

\begin{proof}
We may fix a dense subset $D_y$ of $f^{-1}(y)$ with $1 \le |D_y| \le \lambda$ for each $y \in E$.
We claim that $F = \bigcup \{D_y : y \in E\}$ is as required.

Indeed, then $|F| = \lambda$ and $F \in \mathcal{D}(f^{-1}[E])$ are obvious.
Also, if $H \subs F$ and $|H| < \lambda$ then $f[H] \subs E$
has cardinality $< \lambda$ as well, hence $f[H]$ is not dense in $E$, consequently $H$ cannot
be dense in $F$ by the continuity of $f$. But this means that $d(F) = \lambda$.
\end{proof}

Note that if $f^{-1}[E]$ happens to be dense in $X$ in the above Lemma~\ref{lm:EF},
then so is $F$, hence we have $\lambda \in dd(X)$.

Now, if $f : X \to Y$ is a continuous surjection with $Y$ a $T_1$-space and $y \in Y$ has pseudo-character $\le \lambda$
in $Y$ then $f^{-1}(y)$ is a closed $G_\lambda$-set in~$X$. Thus, in view of Lemma \ref{lm:EF},
we are naturally led to the following definition.

\begin{definition}\label{df:Gl}
For any infinite cardinal $\lambda$ we denote by $\Gamma(\lambda)$ the class of all spaces
in which every closed $G_\lambda$-set has density $\le \lambda$.
\end{definition}

Note that $X \in \Gamma(\lambda)$ clearly implies $d(X) \le \lambda$.
Also, the class  $\Gamma(\lambda)$ is trivially closed under continuous images. This obsevation,
together with our next result, provides a large vaiety of compacta that belong to $\Gamma(\lambda)$.

\begin{theorem}\label{tm:pr}
If $X$ is the product of at most $2^\lambda$ many compacta of weight $\le \lambda$ then $X \in \Gamma(\lambda)$.
\end{theorem}

\begin{proof}
Assume that $X = \prod \{X_i : i \in I\}$, where $|I | \le 2^\lambda$ and $\mathcal{B}_i$ is a base
of the compact space $X_i$ with $|\mathcal{B}_i| \le \lambda$ for every $i \in I$. Then the family
$\mathcal{B}$ consisting of all elementary open sets of the form
$$\bigcap \{p_i^{-1}[B_i] : i \in a\}$$ with $a \in [I]^{< \omega}$ is a base of $X$. Here $p_i : X\to X_i$ is the projection of~$X$ onto its $i$-th factor space $X_i$.
Note that every member of $\mathcal{B}$ depends only on a finite number of coordinates.

Now, consider any closed $G_\lambda$-set $H$ in the product $X$, so $H = \bigcap \mathcal{U}$ for
a family $\mathcal{U}$ of open subsets of $X$ with $|\mathcal{U}| \le \lambda$.
Since $H$ is compact, for every $U \in \mathcal{U}$ there is a finite subset $\mathcal{B}_U$ of
$\mathcal{B}$ such that
$$
    H \subs \bigcup \mathcal{B}_U = V_U \subs U.
$$
Clearly, then every such $V_U$
depends only on a finite number of coordinates, say $a_U$. Then $H$ depends only on the coordinates in
$J = \bigcup \{a_U : U \in \mathcal{U}\}$, where $|J| \le \lambda$.

Let us denote the projection of $H$ to the subproduct $X_J = \prod \{X_i : i \in J\}$ by $H_J$.
Then we have $$d(H_J) \le w(H_J) \le w(X_J) \le \lambda.$$ Also, it is well-known that $|I \setm J| \le 2^\lambda$
implies $d(X_{I \setm J}) \le \lambda$ as well. But $H$ is clearly homeomorphic to $H_J \times X_{I \setm J}$,
consequently we have $d(H) \le \lambda$, and this completes our proof.
\end{proof}

Let us recall now that a (compact) space is called $\mu$-adic for some infinite cardinal number $\mu$
if it is the continuous image of some power of $A(\mu)$, the Alexandrov one-point compactification
of the discrete space of size $\mu$. The space $X$ is {\em polyadic} if it is
$\mu$-adic for some $\mu \ge \omega$. If $X$ is polyadic then we denote by $\mu(X)$ the smallest $\mu$
such that $X$ is $\mu$-adic. It is easy to see that $\mu(X) \le d(X)$.
Clearly, the class of polyadic spaces is the smallest class of spaces that contains all the $A(\mu)$'s
and is closed under products and continuous images.

We note that $\omega$-adic spaces are the same as dyadic spaces, i.e. the continuous images of the Cantor cubes.

Polyadic spaces were first studied by Mr\'owka~\cite{Mrowka70} in 1970 as a generalisation of dyadic spaces.

\begin{corollary}\label{co:po}
If $X$ is a polyadic space of weight $\kappa$ then $X \in \Gamma(\lambda)$ whenever $d(X) \le \lambda \le \kappa$.
\end{corollary}

\begin{proof}
It is standard and easy to prove that $X$  is actually the continuous image of $A(\mu(X))^\kappa$.
But then Theorem \ref{tm:pr} can be applied to the product $A(\mu(X))^\kappa$ because
$w(A(\mu(X))) = \mu(X) \le d(X) \le \lambda$ and $\kappa = w(X) \le 2^{d(X)} \le 2^\lambda$.
Thus we have $A(\mu(X))^\kappa \in \Gamma(\lambda)$, hence $X \in \Gamma(\lambda)$ as well.
\end{proof}

The next theorem is the main result of this section. It involves the cardinal function $\widehat{c}(X)$, the
``hat version" of cellurarity defined as the smallest cardinal $\kappa$ for which there is no
cellular (i.e. disjoint) family of $\kappa$ open sets in $X$.

Another important ingredient of the proof is the fact that a
continuous surjection $f : X \to Y$ is {\em quasi-open}, i.e. $\Int\, f[U] \ne \emptyset$ for every non-empty open
subset of $U$ of $X$, iff for every $E \in \mathcal{D}(Y)$ we have $f^{-1}[E] \in \mathcal{D}(X)$. In fact,
if $X$ is also compact, hence the map $f$ is also closed, then this already follows from $f$ being {\em pseudo-open},
i.e. such that the inverse image of any dense open subset of $Y$ is dense in $X$, see~1.2 and 1.3 of \cite{JSS}.

\begin{theorem}\label{tm:Gl}
Let $X$ be a compact space and $\lambda$ a {\em regular} cardinal such that $X \in \Gamma(\lambda)$, moreover
$$\lambda = \lambda^{<\, \widehat{c}(X)} < \pi(X).$$
Then $\lambda \in dd(X)$.
\end{theorem}

\begin{proof}
We may assume that $X$ is a subspace of the Tychonov cube $[0,1]^\kappa$, where $\kappa = w(X)$ and will
show that there is $J \subs \kappa$ with $|J| = \lambda$ such that the projection map
$p_J : [0,1]^\kappa \to [0,1]^J$ restricted to $X$ is pseudo-open, and hence qusai-open.
To this end, we first take an elementary submodel $M$ of $H(\vartheta)$ for a large enough regular cardinal $\vartheta$ such that
$\,X,\, \kappa \in M$, $\,|M| = \lambda$, and $M$ is $< \widehat{c}(X)$-closed, i.e.
$[M]^{\mu}\subs M$ for all cardinals ${\mu}< \widehat{c}(X)$.
This is possible because $\lambda = \lambda^{<\, \widehat{c}(X)}$ and $\widehat c(X)$ is regular.
Using the regularity of $\lambda$ we may also choose $M$ to be $< \lambda$-covering. that is such that
$M \cap [M]^{< \lambda}$ is cofinal in $[M]^{< \lambda}$. This can be achieved by building up $M$
as the union of an elementary chain of length $\lambda$ of elementary submodels of size $< \lambda$.

After having $M$, we put  $J=M \cap {\kappa}$
and claim that $p_J\restriction X : X \to Y = p_J[X]$ is pseudo-open, i.e.
if $G$ is dense open in $Y$
then $X \cap p_J^{-1}[G] = (p_J\restriction X)^{-1}[G]$ is dense in $X$.

To see this, let us fix an open base $\mc B \in M$ of $[0,1]$ and let $\mc E$ be the
the family of functions with domain a finite subset of $\kappa$ and range included in $\mc B$.
Clearly, $\mc E \in M$ and $M \vDash \varepsilon \in \mc E$ iff $\varepsilon \in M \cap \mc E$
iff $\varepsilon \in \mc E$ and $\dom(\varepsilon) \subs J$. For every $\varepsilon \in \mc E$
 we denote by $[\varepsilon]$ the elementary open set in $[0,1]^\kappa$
determined by $\varepsilon$. Also, for any $\varepsilon \in M \cap \mc E$ we use $[\varepsilon]_J$ to denote
the projection $p_J \big[[\varepsilon] \big]$, which is an elementary open set in $[0,1]^J$.

Since $\widehat{c}(Y) \le \widehat{c}(X)$, we can chose a collection $\mc F \subs M \cap \mc E$ with $|\mc F| < \widehat{c}(X)$
such that for each $\varepsilon \in \mc F$ we have $\emptyset \ne Y \cap [\varepsilon]_J \subs G$, moreover the family
$\{Y \cap [\varepsilon]_J : \varepsilon \in \mc F\}$ is disjoint and its union is dense in $G$ and hence in $Y$.
Note that $|\mc F| < \widehat{c}(X)$ implies $\mc F \in M$.

Now, the fact that $H = \bigcup \{Y \cap [\varepsilon]_J : \varepsilon \in \mc F\}$ is dense in $Y$ can be reformulated
as follows: For every $\varepsilon \in M \cap \mc E$, if $Y \cap [\varepsilon]_J \ne \emptyset$ then there is some
$\eta \in \mc F$ such that $Y \cap [\varepsilon]_J \cap [\eta]_J \ne \emptyset$. But $Y \cap [\varepsilon]_J \cap [\eta]_J \ne \emptyset$
is clearly equivalent with $X \cap [\varepsilon] \cap [\eta] \ne \emptyset$, consequently the following
staement is satisfied in $M$:
$$\forall\, \varepsilon \in \mc E\,\exists\,\eta \in \mc F\,\big( X \cap [\varepsilon] \ne \emptyset\,\Rightarrow\, X \cap [\varepsilon] \cap [\eta] \ne \emptyset \big).$$
Since all the three parameters of this formula, namely $\mc E,\, \mc F$, and $X$ belong to $M$, by elementarity it is actually true.
But this just means that $X \cap p_J^{-1}[H]$ is dense in $X$, hence so is the larger set $X \cap p_J^{-1}[G] $.

Next we show that $\pi(Y) = \lambda$. To see this, using that $M \cap [M]^{< \lambda}$ is cofinal in $[M]^{< \lambda}$,
it suffices to show that $\{Y \cap [\eta]_J : \eta \in \mc F\}$ is not a $\pi$-base of $Y$ for every $\mc F \in M \cap [\mc E]^{< \lambda}$.
But we know that $|\mc F| < \lambda < \pi(X)$, hence there is $\varepsilon \in \mc E$ that witnesses this, i.e.
$X \cap [\varepsilon] \setm [\eta] \ne \emptyset$ for all $\eta \in \mc F$. By $\mc F \in M$ and elementarity, however, we have such an $\varepsilon \in M \cap \mc E$,
i.e. with $\dom(\varepsilon) \subs J$ as well. This then shows that $\{Y \cap [\eta]_J : \eta \in \mc F\}$ is not a $\pi$-base of $Y$.

Now, $Y$ is compact, hence $\lambda = \pi(Y) \in dd(Y)$ by \cite{JSh}, consequently there is $E \in \mc D (Y)$ such that
$|E| = d(E) = \lambda$. Also, we have $w(Y) = \lambda$, hence for every $y \in E$, in fact for every $y \in Y$,
$f^{-1}(y)$ is a closed $G_\lambda$-set in $X$, where $f = p_J\restriction X : X \to Y$. So, as $X \in \Gamma(\lambda)$,
we have $d(f^{-1}(y)) \le \lambda$ as well. But $f$ is quasi-open, consequently $f^{-1}[E]$ is dense in $X$,
consequently by Lemma \ref{lm:EF} we get a dense subset $F$ of $X$ with $|F| = d(F) = \lambda$, hence $\lambda \in dd(X)$.
\end{proof}

Let $X$ be a compact CCC space such that $d(X) < \mathfrak{c} < \pi(X)$ and $X \in \Gamma(\mathfrak{c})$.
By Theorem \ref{df:Gl} then $\mathfrak{c} \in dd(X)$, provided that $\mathfrak{c}$ is regular.
We do not know if the same conclusion is valid if $\mathfrak{c}$ is singular.

Clearly, there is no such problem if $\mathfrak{c}$ is replaced by $\mathfrak{c}^{+n}$, the $n$th successor of $\mathfrak{c}$,
for $0 < n < \omega$. Since, by Theorem 2.1 of \cite{JvMSSz} $dd(X)$ is $\omega$-closed, it follows that the above
observation is also valid if we replace $\mathfrak{c}$ with the singular $\mathfrak{c}^{+\omega}$.

Now, how about $\mathfrak{c}^{+\omega+1}$? As is known, the {\em singular cardinal hypothesis} implies that $\mathfrak{c}^{+\omega+1}$
is an $\omega$-power, hence Theorem \ref{tm:Gl} applies to it. More generally, this argument yields that if $X$ is a compact CCC space
and $d(X) < \mathfrak{c}^{+\alpha} < \pi(X)$ with $X \in \Gamma((\mathfrak{c})^{+\alpha})$ then $\mathfrak{c}^{+\alpha} \in dd(X)$
for all $0 < \alpha < \omega_1$, provided that the singular cardinal hypothesis holds. It is well-known that the failure of the
singular cardinal hypothesis needs large cardinals. We do not know, however, if the failure of the above statement say for
$\mathfrak{c}^{+\omega + 1}$, or for $\mathfrak{c}^{+\alpha}$ with $0 < \alpha < \omega_1$, would imply the existence of
some large cardinal.

Finally we note that any dyadic compactum $X$ is CCC and belongs to $\Gamma(\lambda)$ for any $d(X) \le \lambda \le w(X) = \pi(X)$
by Corollary \ref{co:po}. In the next section we shall show that in this case we always have $\lambda \in dd(X)$,
i.e. $dd(X) = [d(X), w(X)]$.

\section{The class $\Delta(\lambda)$ and $dd(X)$ for polyadic $X$}

In this section for any infinite cardinal $\lambda$ we introduce a class of {\em compact} spaces that will turn out
to include all polyadic spaces $X$ such that $d(X) \le \lambda \le w(X)$. The significance of these classes will become clear
later.

\begin{definition}\label{df:Dl}
For any infinite cardinal $\lambda$ we denote by $\Delta(\lambda)$ the class of all {\em compact} spaces of weight $\ge \lambda$
in which there is a dense subset $S$ such that \\ (i) $w(\overline{T}) < \lambda$ for every $T \subs S$ with $|T| < \lambda$;\\
(ii) if $f : X \to Y$ is continuous onto and $w(Y) = \lambda$ then $|f[S]| = \lambda$.
\end{definition}

While this definition might seem ad hoc, our following result justifies its relevance. We also note here that every class $\Delta(\lambda)$
is closed under continuous images of weight $\ge \lambda$. Indeed, assume that $X \in \Delta(\lambda)$ and~$S$ is a dense subset of $X$ witnessing this.
It is straight forward to see then that if $f : X \to Y$ is a continuous surjection
with $w(Y) \ge \lambda$, then $f[S]$ is a witness for $Y \in \Delta(\lambda)$.

\begin{theorem}\label{tm:GD}
If $X \in \Gamma(\lambda) \cap \Delta(\lambda)$ then $\lambda \in dd(X)$.
\end{theorem}

\begin{proof}
Let $S$ be a dense subset of $X$ witnessing $X \in \Delta(\lambda)$. By \cite{Juhasz97} there is a continuous
surjection $f : X \to Y$ such that $w(Y) = \lambda$ because $w(X) \ge \lambda$, and set $E = f[S]$.
Then $E$ is dense in $Y$ and by clause (ii) of Definition \ref{df:Gl} we have $|E| = \lambda$.

On the other hand, by (i) we have $w(f[\overline{T}]) \le w(\overline{T}) < \lambda = w(Y)$ for any $T \in [S]^{< \lambda}$,
because the weight of a continuous image of a compact space may not be bigger than the weight of the domain.
Consequently, no $D \in [E]^{< \lambda}$  is dense in $E$ (or equivalently $Y$) because there is $T \in [S]^{< \lambda}$ with $f[T] = D$,
and so $w(\overline{D}) \le w(\overline{T}) < \lambda = w(Y)$. In other words, we have $|E| = d(E) = \lambda$.

Now, $w(Y) = \lambda$ clearly implies that $f^{-1}(y)$ is a closed $G_\lambda$-set in $X$ for all $y \in Y$, hence
$X \in \Gamma(\lambda)$ implies $d(f^{-1}(y)) \le \lambda$.
Thus we may apply Lemma \ref{lm:EF} to conclude that there is a dense subset
$F$ of $f^{-1}[E]$ such that $|F| = d(F) = \lambda$. But we clearly have $S \subs f^{-1}[E]$, hence $f^{-1}[E]$ is dense in $X$,
consequently so is $F$. This , however, means that $F$ is a witness for $\lambda \in dd(X)$.
\end{proof}

Our aim now is to show that all polyadic spaces belong to $\Delta(\lambda)$ for the relevant values of $\lambda$.
To do that, we need some preparation.

We start with the trivial observation that any space $X$ decomposes in the disjoint union $X = I(X) \cup X'$,
where $I(X)$ is the set of all isolated points and $X'$ is the setof all accumulation points of $X$.
The following simple observation will turn out to be crucial below.

\begin{lemma}\label{lm:IX}
Let $f : X \to Y$ be any closed and continuous map. Then $f[X] \setm f[X'] \subs I(f[X])$, and so
it is a discrete subset of $Y$.
\end{lemma}

\begin{proof}
Indeed, let $y \in f[X] \setm f[X']$, this means that $f^{-1}(y) \subs I(X)$. Consequently,
$X \setm f^{-1}(y)$ is closed in $X$, hence its $f$-image $f[X] \setm \{y\}$ is closed in $f[X]$,
hence $y \in I(f[X])$.
\end{proof}

The next result is the most significant step towards our aim.

\begin{theorem}\label{tm:mk}
If $\omega \le \mu < \lambda \le \kappa$ then $A(\mu)^\kappa \in \Delta(\lambda)$.
\end{theorem}

\begin{proof}
To make our notation specific, we assume that the underlying set of $A(\mu)$ is $\mu + 1 = \mu \cup \{\mu\}$,
with all $\alpha < \mu$ isolated and the top ordinal $\mu$ being the compactifying point.
Then the points of $A(\mu)^\kappa$ are functions with domain $\kappa$ and range included in $\mu + 1$.
Clearly, $A(\mu)^\kappa$ is compact of weight $\ge \lambda$, hence we only have to find a dense $S \subs A(\mu)^\kappa$
satisfying clauses (i) and (ii) of Definition \ref{df:Dl}.

To do that, for every $x \in A(\mu)^\kappa$ we define the support $\supp(x)$ of $x$ as follows:
$$
    \supp(x) = \{\xi < \kappa : x(\xi) \in \mu\}.
$$
We then set $S = \{x \in A(\mu)^\kappa : |\supp(x) < \omega|\}$,
clearly $S$ is dense in $A(\mu)^\kappa$ and we claim that it witnesses $A(\mu)^\kappa \in \Delta(\lambda)$.
Interestingly, $S$ does not depend on $\lambda$.

To check (i), take any $T \in [S]^{< \lambda}$. Then $\nu = |\bigcup \{\supp(x) : x \in T\}| < \lambda$ as well,
hence $\overline{T}$ is homeomorphic to a subspace of $A(\mu)^\nu$ that clearly has weight $\mu \cdot \nu < \lambda$.

Now we turn to checking (ii). To this end, we let $S_n = \{x \in S : |\supp(x) \le n|\}$ for any $n < \omega$.
Every such $S_n$ is clearly closed in $A(\mu)^\kappa$, and hence is compact. Also, it is clear that
$$
    I(S_{n+1}) = S_{n+1} \setm S_n = \{x \in S : |\supp(x)| = n + 1\}.
$$
Note that $S_0$ is the singleton consisting of
the constant $\mu$-valued function on $\kappa$.

Let now $f$ be a continuous surjection of $A(\mu)^\kappa$ onto $Y$ with $w(Y) = \lambda$. We shall show by induction on
$n < \omega$ that $|f[S_n]| \le \lambda$ for all $n < \omega$. This is trivial for $n = 0$. The induction step from
$n$ to $n + 1$ uses Lemma \ref{lm:IX} for the obviously closed map $f \upharpoonright S_{n +1}$.
Indeed, this tells us that $f[S_{n+1})] \setm f[S_n]$ is discrete in $Y$, and hence has cardinality $\le \lambda$,
hence $|f[S_n]| \le \lambda$ implies $|f[S_{n+1}]| \le \lambda$. This, together with $S = \bigcup \{S_n : n < \omega\}$
then implies  $|f[S]| \le \lambda$ as well. But (i) clearly implies that $f[T]$ cannot be dense in $Y$ for $T \in [S]^{< \lambda}$,
hence we actually have $|f[S]| = \lambda$, completing the proof of (ii).
\end{proof}

Now we can easily formulate and prove what we aimed for.

\begin{theorem}\label{tm:pdd}
If $X$ is a polyadic space then $dd(X) = [d(X), w(X)]$.
\end{theorem}

\begin{proof}
If $d(X) = w(X)$ then the statement is obvious, so assume that $d(X) < w(X)$. As we have seen in the previous section,
$\mu(X) \le d(X)$ and $X$ is the continuous image of $A(\mu(X))^{w(X)}$. Since we always have $d(X) = \min dd(X)$, it remains
to consider $\lambda$ with $\mu(X) \le d(X) < \lambda \le w(X)$. By Theorem \ref{tm:mk} in this case $A(\mu(X))^{w(X)} \in \Delta(\lambda)$,
hence $\lambda \le w(X)$ implies $X \in \Delta(\lambda)$ as well. However, by Corollary \ref{co:po} we also have $X \in \Gamma(\lambda)$,
consequently Theorem \ref{tm:GD} implies that $\lambda \in dd(X)$.
\end{proof}

It is worth to note that, as a byproduct, this result implies that $\pi(X) = w(X)$ for any polyadic space $X$ because
$\pi(X) = \max dd(X) \le w(X)$. Of course, this fact is well-known.

As is well-known, any compact topological group $G$ is dyadic, consequently we have $dd(G) = [d(G), w(G)] = [\log \,w(G), w(G)]$.
In fact, the equality $dd(G) = [d(G), w(G)]$ holds for all locally compact topological groups as well, being an immediate
consequence of the following result.

\begin{theorem}\label{tm:lcg}
The Alexandrov one-point compactification $A(G)$ of any locally compact topological group $G$
is polyadic.
\end{theorem}

\begin{proof}
By Clearly and Morris~\cite{ClearyMorris88}, we may assume that $$G = \mathbb{R}^n \times K\times D,$$ where $n<\omega,$ $K$ is a compact group, and $D$ is discrete.
It suffices to prove that $G$ has a polyadic compactification $aG$ since the 1-point compactification of $G$ is a continuous image of $aG$.
The 1-point compactification $A(\mathbb{R}^n)$ of $\mathbb{R}^n$ is dyadic, being compact metric.
The 1-point compactification $A(D)$ is trivially polyadic since $D$ is discrete. Hence as $K$ is dyadic,
and the product of polyadic spaces is polyadic, the compactification
$A(\mathbb{R}^n) \times K \times A(D)$ of $G$ is indeed polyadic.
\end{proof}

Since it is obvious that $dd(G) = dd(A(G))$ and $\pi(G) = \pi(A(G)) = w(A(G))$,
we indeed conclude that $dd(G) = [d(G), \pi(G)] = [d(G), w(G)]$.


\begin{thebibliography}{12}


\bibitem{ClearyMorris88}
J.~Cleary and S.~Morris, {\em Topologies on locally compact groups}, Bull.
  Austr. Math. Soc. {\bf 38} (1988), 105--111.

\bibitem{juhasz}
I.~Juh\'asz, {\em {Cardinal functions in topology}}, Mathematical Centre Tract,
  vol.~34, Mathematical Centre, Amsterdam, 1971.

\bibitem{Juhasz97}
I.~Juh\'{a}sz, {\em Cardinal functions on continuous images}, Proceedings of
  the {T}ennessee {T}opology {C}onference ({N}ashville, {TN}, 1996), World Sci.
  Publ., River Edge, NJ, pp.~89--94.

\bibitem{JvMSSz}
I.~Juh\'asz, J.~van Mill, L.~Soukup, and Z.~Szentmikl\'ossy,
The double density spectrum of a topological space, Israel J. Math.,\\
published online in December 2022: https://rdcu.be/c2ml4

\bibitem{JSh}
I. Juhász and S. Shelah, {\em {$\pi(X)=\delta(X)$ for compact $X$}}, Top.
  Appl. {\bf 32} (1989), 289--294.

\bibitem{JSS}
I.~Juh{\'a}sz, L.~Soukup, and Z.~Szentmikl{\'o}ssy,
{\em Spaces of small cellularity have nowhere constant continuous images  of small weight},
Top. Appl., {\bf 281} (2020), 107212

\bibitem{Mrowka70}
S. Mr\'owka, {\em Mazur theorem and $m$-adic spaces}, Bull. Polon. Acad. Sci. S\'{e}r. Math. Astronom. Phys., {\bf 18} (1970), 299--305.

\end{thebibliography}
\end{document}